\documentclass[12pt]{amsart}
\usepackage{amssymb}
\usepackage[all]{xy}
\newtheorem{theorem}{Theorem}[]

\newtheorem{proposition}[theorem]{Proposition}
\newtheorem{definition}[theorem]{Definition}
\newtheorem{question}[theorem]{Question}
\newtheorem{remark}[theorem]{Remark}
\def\P{\Phi}
\def\S{\Psi}

\def\G{\Gamma}
\def\D{\Delta}
\def\i{\mathrm{I}}
\def\Q{\mathfrak{Q}}
\def\to{\longrightarrow}
\def\ot{\otimes}

\begin{document}
%%%%%%%%%%%%%%%%%%%%%%%%%%%%%%%%%%%%%%%%%%%%%%%%%%%%%%%%%%%%%%%%%%%%%%%%%%%%%%%%%%%%%%%%%%%%%%%%%%%%%%%%%%%%%%%%%%%%%%%%%%
\title[Compact Quantum Semigroups]{A Kind of Compact Quantum Semigroups}
%%%%%%%%%%%%%%%%%%%%%%%%%%%%%%%%%%%%%%%%%%%%%%%%%%%%%%%%%%%%%%%%%%%%%%%%%%%%%%%%%%%%%%%%%%%%%%%%%%%%%%%%%%%%%%%%%%%%%%%%%%
\author[M. M. Sadr]{Maysam Maysami Sadr}
\address{Depertment of Mathematics\\
Institute for Advanced Studies in Basic Sciences (IASBS of Zanjan) \\
P.O. Box 45195-1159, Zanjan 45137-66731, Iran}
\email{sadr@iasbs.ac.ir}
%%%%%%%%%%%%%%%%%%%%%%%%%%%%%%%%%%%%%%%%%%%%%%%%%%%%%%%%%%%%%%%%%%%%%%%%%%%%%%%%%%%%%%%%%%%%%%%%%%%%%%%%%%%%%%%%%%%%%%%%%%
\begin{abstract}
We show that the quantum family of all maps from a finite space to a finite dimensional
compact quantum semigroup has a canonical quantum semigroup structure.
\end{abstract}
%%%%%%%%%%%%%%%%%%%%%%%%%%%%%%%%%%%%%%%%%%%%%%%%%%%%%%%%%%%%%%%%%%%%%%%%%%%%%%%%%%%%%%%%%%%%%%%%%%%%%%%%%%%%%%%%%%%%%%%%%%
\maketitle
%%%%%%%%%%%%%%%%%%%%%%%%%%%%%%%%%%%%%%%%%%%%%%%%%%%%%%%%%%%%%%%%%%%%%%%%%%%%%%%%%%%%%%%%%%%%%%%%%%%%%%%%%%%%%%%%%%%%%%%%%%
\section{Introduction}
According to the Gelfand duality, the category of compact Hausdorff spaces and
continuous maps and the category of commutative unital C*-algebras and unital *-homomorphisms are dual.
In this duality, any compact space $X$ corresponds to $\mathcal{C}(X)$, the C*-algebra of all continuous complex valued maps on $X$,
and any commutative unital C*-algebra corresponds to its maximal ideal space. Thus as the fundamental concept in
{\it Non-Commutative Topology}, a non-commutative unital C*-algebra $A$ is considered as the algebra of
continuous functions on a {\it symbolic} compact non-commutative space $\Q A$.
In this correspondence, *-homomorphisms $\P:A\to B$ interpret as symbolic
continuous maps $\Q\P:\Q B\to\Q A$. Since the coordinates observable of a quantum mechanical systems are non-commutative, 
some times  non-commutative spaces are called quantum spaces.  

Woronowicz \cite{W1} and So{\l}tan \cite{S1} have defined a quantum space $\Q C$ of {\it all maps}
from $\Q B$ to $\Q A$ and showed that $C$ exists under appropriate conditions on $A$ and $B$.
In \cite{S}, we considered the functorial properties of this notion.
In this short note, we show that if $\Q A$ is a compact finite dimensional
(i.e. $A$ is unital and finitely generated) quantum semigroup,
and if $\Q B$ is a finite commutative quantum space (i.e. $B$ is a finite dimensional commutative
C*-algebra), then $\Q C$ has a canonical quantum semigroup structure. In the other words, we construct
the non-commutative version of semigroup $\mathcal{F}(X,S)$ described as follows:

{\it Let $X$ be a finite space and $S$ be a compact semigroup. Then the space $\mathcal{F}(X,S)$
of all maps from $X$ to $S$, is a compact semigroup with compact-open topology and pointwise multiplication.}

In Section 2, we define quantum families of all maps and compact quantum semigroups.
In  Section 3, we state and prove our main result; also we consider a result about quantum semigroups with counits.
At least, in Section 4, we consider some examples.
%%%%%%%%%%%%%%%%%%%%%%%%%%%%%%%%%%%%%%%%%%%%%%%%%%%%%%%%%%%%%%%%%%%%%%%%%%%%%%%%%%%%%%%%%%%%%%%%%%%%%%%%%%%%%%%%%%%%%%%%%%
\section{Quantum families of maps and quantum semigroups}
All C*-algebras in this paper have unit and all C*-algebra homomorphisms preserve the units.
For any C*-algebra $A$, $\i_A$ and $1_A$ denote the identity homomorphism from $A$ to $A$, and the unit of $A$, respectively.
For C*-algebras $A,B$, $A\otimes B$ denotes the spatial tensor product of $A$ and $B$.
If $\P:A\to B$ and $\P':A'\to B'$ are *-homomorphisms, then $\P{\ot}\P':A\otimes A'\to B\otimes B'$ is the
*-homomorphism defined by $\P{\ot}\P'(a{\ot}a')=\P(a){\ot}\P'(a')$ ($a\in A,a'\in A'$).

Let $X,Y$ and $Z$ be three compact Hausdorff spaces and $\mathcal{C}(Y,X)$ be the space of all continuous maps from $Y$ to $X$ with compact open topology. Consider a continuous map $f:Z\to\mathcal{C}(Y,X)$. Then the pair $(Z,f)$ is
a continuous {\it family of maps} from $Y$ to $X$ indexed by $f$ with parameters in $Z$. On the other hand, by topological {\it exponential law} we know that $f$ is characterized by a continuous map $\tilde{f}:Y\times Z\to X$
defined by $\tilde{f}(y,z)=f(z)(y)$. Thus $(Z,\tilde{f})$ can be considered as a family of maps from $Y$ to $X$.
Now, by Gelfand's duality we can simply translate this system to non-commutative language:
\begin{definition}\label{def1}
(\cite{W1},\cite{S1})
Let $A$ and $B$ be unital C*-algebras. By a quantum family of maps from $\Q B$ to $\Q A$, we mean a pair $(C,\P)$,
containing a unital C*-algebra $C$ and a unital *-homomorphism $\P:A\to B{\otimes} C$.
\end{definition}
Now, suppose  instead of parameter space $Z$ we use $\mathcal{C}(Y,X)$ (note that in general this space is
not locally compact). Then the family
$$\i\mathrm{d}:\mathcal{C}(Y,X)\to\mathcal{C}(Y,X)\hspace{3mm}(\tilde{\i\mathrm{d}}:\mathcal{C}(Y,X)\times Y\to X)$$
of {\it all} maps from $Y$ to $X$ has the following universal property:\\
For every family $\tilde{f}:Z\times Y\to X$ of maps from $Y$ to $X$, there is a unique map $f:Z\to\mathcal{C}(Y,X)$
such that the following diagram is commutative.
\[\xymatrix{Z\times Y\ar[rr]^-{\tilde{f}}\ar[d]^{f\times\i\mathrm{d}_Y}&& X\ar@{=}[d]\\
\mathcal{C}(Y,X)\times Y\ar[rr]^-{\tilde{\i\mathrm{d}}}&& X}\]
Thus, we can make the following Definition in non-commutative setting:
\begin{definition}\label{def2}
(\cite{W1},\cite{S1})
With the assumptions of  Definition \ref{def1}, $(C,\P)$ is called a quantum family of all maps from $\Q B$ to $\Q A$
if for every unital C*-algebra $D$ and any unital *-homomorphism $\S:A\to B{\otimes} D$, there is a unique unital
*-homomorphism $\G:C\to D$ such that the following diagram is commutative.
\[\xymatrix{A\ar[rr]^-{\P}\ar@{=}[d]&& B\ot C\ar[d]^{\i_B{\ot}\G}\\
A\ar[rr]^-{\S}&& B\ot D}\]
\end{definition}
By the universal property of Definition \ref{def2}, it is clear that if $(C,\P)$ and $(C',\P')$
are two quantum families of all maps from $\Q B$ to $\Q A$, then there is a *-isometric isomorphism
between $C$ and $C'$.
\begin{proposition}
Let $A$ be a unital finitely generated C*-algebra and $B$ be a finite dimensional C*-algebra. Then
the quantum family of all maps from $\Q B$ to $\Q A$ exists.
\end{proposition}
\begin{proof}
See \cite{W1} or \cite{S1}.
\end{proof}
\begin{definition}
(\cite{W2},\cite{MurphyTuset1},\cite{S1},\cite{S2}) A pair $(A,\D)$ consisting of a unital C*-algebra $A$ and a unital
*-homomorphism $\D:A\to A{\otimes} A$ is called a compact quantum semigroup if $\D$ is a coassociative
comultiplication: $(\D{\ot} \i_A)  \D=(\i_A{\ot}\D)  \D$.
\end{definition}
A *-homomorphism $\D:A\to A{\otimes} A$ induces a binary operation on the dual space $A^*$ defined by
$\tau\sigma=(\tau\ot\sigma)\D$ for $\tau,\sigma\in A^*$.
Now, suppose that $S$ is a compact Hausdorff topological semigroup. Using the canonical identity
$\mathcal{C}(S)\ot\mathcal{C}(S)\cong\mathcal{C}(S\times S)$, we define a *-homomorphism
$\D:\mathcal{C}(S)\to\mathcal{C}(S)\otimes\mathcal{C}(S)$ by
$\D(f)(s,s')=f(ss')$ for $f\in\mathcal{C}(S)$ and $s,s'\in S$.
Then $\D$ is a coassociative comultiplication on $\mathcal{C}(S)$ and thus $(\mathcal{C}(S),\D)$ is a compact quantum
semigroup. Conversely, if $(A,\D)$ is a compact quantum semigroup such that $A$ is abelian, then the character space
of $A$, with the binary operation induced by $\D$, is a compact Hausdorff topological semigroup (\cite{VV}).
It is well known that a compact semigroup with cancelation property is a compact group (\cite[Proposition 3.2]{MaesVandaele1}).
Analogous cancelation properties for quantum semigroups are defined as follows.
\begin{definition}
Let $(A,\D)$ be a compact quantum semigroup.
\begin{enumerate}
\item [(i)] (\cite{MurphyTuset1}) $(A,\D)$ has left (resp. right) cancelation property if the linear span of $\{(b\ot1)\D(a):a,b\in A\}$
(resp. $\{(1\ot b)\D(a):a,b\in A\}$) is dense in $A\ot A$.
\item [(ii)] (\cite{MurphyTuset1}) $(A,\D)$ has weak left cancelation property if, whenever $\tau,\sigma\in A^*$ are such that
$(\tau a)\sigma=0$ for all $a\in A$, we must have $\tau=0$ or $\sigma=0$. Similarly, $(A,\D)$ has weak right cancelation property if,
whenever $\tau(\sigma a)=0$ for all $a\in A$, we must have $\tau=0$ or $\sigma=0$.
\item [(iii)] (\cite{S1}) A left (resp. right) counit for $(A,\D)$, is a character $\epsilon$ on $A$ (a unital *-homomorphism
$\epsilon:A\to\mathbb{C}$), satisfying $(\epsilon\ot\i_A)\D=\i_A$ (resp. $(\i_A\ot\epsilon)\D=\i_A$. A left and right counit is called
(two-sided) counit.
\end{enumerate}
\end{definition}
In the above Definition the functionals $\tau a$ and $a\tau$ are defined by $\tau a(x)=\tau(ax)$ and $a\tau(x)=\tau(xa)$.
\begin{remark}
In \cite{W2}, counits are characters on special dense subalgebras of compact quantum groups. These subalgebras are constructed from
finite dimensional unitary representations of compact quantum groups. In this paper we mainly deal with quantum semigroups and since it is not
natural to define unitary representations for (quantum) semigroups, we use the above notion for counits.
\end{remark}
It is clear that the left (resp. right) cancelation property implies weak left (resp. weak right) cancelation property.
The converse is partially satisfied (\cite[Theorem 3.2]{MurphyTuset1}):
\begin{theorem}
Let $(A,\D)$ be a compact quantum semigroup. Then $(A,\D)$ has both left and right cancelation properties if and only if
it has both weak left and weak right cancelation properties.
\end{theorem}
\begin{definition}
(\cite{W2},\cite{MurphyTuset1},\cite{MaesVandaele1}) A compact quantum semigroup with both left and right cancelation properties
is called compact quantum group.
\end{definition}
Again consider compact semigroup $S$ and its corresponding compact quantum semigroup $(\mathcal{C}(S),\D)$
defined above. Using Proposition 3.2 of \cite{MaesVandaele1}, it is easily proved that $S$ is a compact group
if and only if $(\mathcal{C}(S),\D)$ is a compact quantum group.
%%%%%%%%%%%%%%%%%%%%%%%%%%%%%%%%%%%%%%%%%%%%%%%%%%%%%%%%%%%%%%%%%%%%%%%%%%%%%%%%%%%%%%%%%%%%%%%%%%%%%%%%%%%%%%%%%%%%%%%%%%
\section{The results}
In this section, we state and prove the main result.
\begin{theorem}\label{t1}
Let $(A,\D)$ be a compact quantum semigroup with finitely generated $A$, $B$ be a finite dimensional
commutative C*-algebra, and $(C,\P)$ be the quantum family of all maps from $\Q B$ to $\Q A$. Consider
the unique unital *-homomorphism $\G:C\to C\otimes C$
such that the diagram
\begin{equation}\label{d1}
\xymatrix{A\ar[rr]^-{\P}\ar[d]^{\D}&& B\ot C\ar[d]^{\i_B{\ot}\G}\\
A\ot A\ar[d]^{\P{\ot}\P}&& B\ot C\ot C\\
B\ot C\ot B\ot C\ar[rr]^-{\i_B\ot F\ot \i_C}&&B\ot B\ot C\ot C\ar[u]^{m{\ot} \i_{C{\ot} C}}}
\end{equation}
is commutative, where $F:C{\otimes} B\to B{\otimes} C$ is the flip map, i.e.
$c{\ot} b\longmapsto b{\ot} c$ ($b\in B, c\in C$), and
$m:B{\otimes} B\to B$ is the multiplication *-homomorphism of $B$, i.e. $m(b{\ot} b')=bb'$ ($b,b'\in B)$.
Then $(C,\G)$ is a compact quantum semigroup.
\end{theorem}
\begin{proof}
We must prove that $(\i_C{\ot}\G) \G=(\G{\ot}\i_C) \G$, and for this, by the universal property of quantum
families of maps, it is enough to prove that
\begin{equation}\label{e1}
(\i_B{\ot}\i_C{\ot}\G) (\i_B{\ot}\G) \P=(\i_B{\ot}\G{\ot}\i_C) (\i_B{\ot}\G) \P.
\end{equation}
Note that by the commutativity of (\ref{d1}), we have
$$(\i_B{\ot}\G) \P=(m{\ot} \i_{C{\ot} C}) (\i_B{\ot} F{\ot} \i_C) (\P{\ot}\P) \D.$$
Let us begin from the left hand side of (\ref{e1}):
\begin{equation*}
\begin{split}
&(\i_B{\ot}\i_C{\ot}\G) (\i_B{\ot}\G) \P\\
=&(\i_B{\ot}\i_C{\ot}\G) (m{\ot} \i_{C{\ot} C}) (\i_B{\ot} F{\ot} \i_C) (\P{\ot}\P) \D\\
=&(m{\ot}\i_C{\ot}\G) (\i_B{\ot} F{\ot} \i_C) (\P{\ot}\P) \D\\
=&(m{\ot}\i_{C{\ot}C{\ot}C})(\i_B{\ot}F{\ot}\i_{C{\ot}C})(\i_{B{\ot}C{\ot}B}{\ot}\G)(\P{\ot}\P)\D\\
=&(m{\ot}\i_{C{\ot}C{\ot}C})(\i_B{\ot}F{\ot}\i_{C{\ot}C})(\i_{B{\ot}C{\ot}B}{\ot}\G)(\P{\ot}\i_{B{\ot}C})
(\i_A{\ot}\P)\D\\
=&(m{\ot}\i_{C{\ot}C{\ot}C})(\i_B{\ot}F{\ot}\i_{C{\ot}C})(\P{\ot}\i_{B{\ot}C{\ot}C})(\i_A{\ot}\i_B{\ot}\G)
(\i_A{\ot}\P)\D\\
=&(m{\ot}\i_{C{\ot}C{\ot}C})(\i_B{\ot}F{\ot}\i_{C{\ot}C})(\P{\ot}\i_B{\ot}\G)(\i_A{\ot}\P)\D\\
=&(m{\ot}\i_{C{\ot}C{\ot}C})(\i_B{\ot}F{\ot}\i_{C{\ot}C})(\P{\ot}[(\i_B{\ot}\G)\P])\D\\
=&(m{\ot}\i_{C{\ot}C{\ot}C})(\i_B{\ot}F{\ot}\i_{C{\ot}C})
(\P{\ot}[(m{\ot} \i_{C{\ot} C}) (\i_B{\ot} F{\ot} \i_C) (\P{\ot}\P) \D])\D\\
=&(m{\ot}\i_{C{\ot}C{\ot}C})(\i_B{\ot}F{\ot}\i_{C{\ot}C})(\i_{B{\ot}C}{\ot}m{\ot} \i_{C{\ot} C})
(\i_{B{\ot}C{\ot}B}{\ot} F{\ot} \i_C)(\P{\ot}\P{\ot}\P)(\i_A{\ot}\D)\D
\end{split}
\end{equation*}
For the right hand side of (\ref{e1}), we have
\begin{equation*}
\begin{split}
&(\i_B{\ot}\G{\ot}\i_C) (\i_B{\ot}\G) \P\\
=&(\i_B{\ot}\G{\ot}\i_C)(m{\ot} \i_{C{\ot} C}) (\i_B{\ot} F{\ot} \i_C) (\P{\ot}\P) \D\\
=&(m{\ot}\G{\ot}\i_C)(\i_B{\ot}F{\ot}\i_C)(\P{\ot}\P) \D\\
=&(m{\ot}\i_{C{\ot}C{\ot}C})(\i_B{\ot}F{\ot}\i_C{\ot}\i_C)(\i_B{\ot}\i_C{\ot}F{\ot}\i_C)
(\i_B{\ot}\G{\ot}\i_B{\ot}\i_C)(\P{\ot}\P) \D\\
=&(m{\ot}\i_{C{\ot}C{\ot}C})(\i_B{\ot}F{\ot}\i_{C{\ot}C})(\i_{B{\ot}C}{\ot}F{\ot}\i_C)([(\i_B{\ot}\G) \P]{\ot}\P)\D\\
=&(m{\ot}\i_{C{\ot}C{\ot}C})(\i_B{\ot}F{\ot}\i_{C{\ot}C})(\i_{B{\ot}C}{\ot}F{\ot}\i_C)
([(m{\ot} \i_{C{\ot} C}) (\i_B{\ot} F{\ot} \i_C) (\P{\ot}\P) \D]{\ot}\P)\D,
\end{split}
\end{equation*}
and thus if $W=(\i_B{\ot}F{\ot}\i_{C{\ot}C})(\i_{B{\ot}C}{\ot}F{\ot}\i_C)$, then
\begin{equation*}
\begin{split}
&(\i_B{\ot}\G{\ot}\i_C) (\i_B{\ot}\G) \P=\\
&(m{\ot}\i_{C{\ot}C{\ot}C})W
(m{\ot}\i_{C{\ot}C{\ot}B{\ot}C})(\i_B{\ot}F{\ot}\i_{C{\ot}B{\ot}C})(\P{\ot}\P{\ot}\P)(\D{\ot}\i_A)\D.
\end{split}
\end{equation*}
Thus, since $(\i_A{\ot}\D)\D=(\D{\ot}\i_A)\D$, to prove (\ref{e1}), it is enough to show that
\begin{equation}\label{e2}
\begin{split}
&(m{\ot}\i_{C{\ot}C{\ot}C})(\i_B{\ot}F{\ot}\i_{C{\ot}C})(\i_{B{\ot}C}{\ot}m{\ot} \i_{C{\ot} C})
(\i_{B{\ot}C{\ot}B}{\ot} F{\ot} \i_C)=\\
&(m{\ot}\i_{C{\ot}C{\ot}C})W(m{\ot}\i_{C{\ot}C{\ot}B{\ot}C})(\i_B{\ot}F{\ot}\i_{C{\ot}B{\ot}C}).
\end{split}
\end{equation}
Let $b_1,b_2,b_3\in B$ and $c_1,c_2,c_3\in C$. Then for the left hand side of (\ref{e2}), we have,
\begin{equation*}
\begin{split}
&(m{\ot}\i_{C{\ot}C{\ot}C})(\i_B{\ot}F{\ot}\i_{C{\ot}C})(\i_{B{\ot}C}{\ot}m{\ot} \i_{C{\ot} C})
(\i_{B{\ot}C{\ot}B}{\ot} F{\ot} \i_C)(b_1{\ot}c_1{\ot}b_2{\ot}c_2{\ot}b_3{\ot}c_3)\\
&=(m{\ot}\i_{C{\ot}C{\ot}C})(\i_B{\ot}F{\ot}\i_{C{\ot}C})(\i_{B{\ot}C}{\ot}m{\ot} \i_{C{\ot} C})
(b_1{\ot}c_1{\ot}b_2{\ot}b_3{\ot}c_2{\ot}c_3)\\
&=(m{\ot}\i_{C{\ot}C{\ot}C})(\i_B{\ot}F{\ot}\i_{C{\ot}C})(b_1{\ot}c_1{\ot}(b_2b_3){\ot}c_2{\ot}c_3)\\
&=(m{\ot}\i_{C{\ot}C{\ot}C})(b_1{\ot}(b_2b_3){\ot}c_1{\ot}c_2{\ot}c_3)\\
&=b_1(b_2b_3){\ot}c_1{\ot}c_2{\ot}c_3\\
&=(b_1b_2b_3){\ot}c_1{\ot}c_2{\ot}c_3,
\end{split}
\end{equation*}
and for the right hand side of (\ref{e2}),
\begin{equation*}
\begin{split}
&(m{\ot}\i_{C{\ot}C{\ot}C})W(m{\ot}\i_{C{\ot}C{\ot}B{\ot}C})(\i_B{\ot}F{\ot}\i_{C{\ot}B{\ot}C})
(b_1{\ot}c_1{\ot}b_2{\ot}c_2{\ot}b_3{\ot}c_3)\\
=&(m{\ot}\i_{C{\ot}C{\ot}C})W(m{\ot}\i_{C{\ot}C{\ot}B{\ot}C})(b_1{\ot}b_2{\ot}c_1{\ot}c_2{\ot}b_3{\ot}c_3)\\
=&(m{\ot}\i_{C{\ot}C{\ot}C})W(b_1b_2{\ot}c_1{\ot}c_2{\ot}b_3{\ot}c_3)\\
=&(m{\ot}\i_{C{\ot}C{\ot}C})(b_1b_2{\ot}b_3{\ot}c_1{\ot}c_2{\ot}c_3)\\
=&(b_1b_2b_3){\ot}c_1{\ot}c_2{\ot}c_3
\end{split}
\end{equation*}
Therefore, (\ref{e2}) is satisfied and the proof is complete.
\end{proof}
\begin{theorem}\label{t2}
Let $(A,\D)$ be a compact quantum semigroup with a left counit. Suppose that $B,C,\P$ and $\G$ are as in Theorem \ref{t1}.
Then the compact quantum semigroup $(C,\G)$ has a left counit.
\end{theorem}
\begin{proof}
Let $\epsilon:A\to\mathbb{C}$ be a left counit for $(A,\D)$. Defin the unital *-algebra homomorphism $\omega:A\to B\ot\mathbb{C}=B$ by
$\omega(a)=1_B\ot\epsilon(a)=\epsilon(a)1_B$ ($a\in A$). Then the universal property of $(C,\P)$ shows that there is a character
$\hat{\epsilon}:C\to\mathbb{C}$ such that the following diagram is commutative:
\[\xymatrix{A\ar[rr]^-{\P}\ar@{=}[d]&& B\ot C\ar[d]^{\i_B{\ot}\hat{\epsilon}}\\
A\ar[rr]^-{\omega}&& B\ot \mathbb{C}}\]
We show that $(\hat{\epsilon}\ot\i_C)\G=\i_C$, and thus $\hat{\epsilon}$ is a counit for $(C,\G)$. By the universal property of $(C,\P)$,
it is enough to show that
\begin{equation}\label{e3}
(\i_B{\ot}[(\hat{\epsilon}\ot\i_C)\G])\P=\P.
\end{equation}
We have
\begin{equation*}
\begin{split}
(\i_B{\ot}[(\hat{\epsilon}\ot\i_C)\G])\P=&(\i_B{\ot}\hat{\epsilon}{\ot}\i_C) (\i_B{\ot}\G)\P\\
=&(\i_B{\ot}\hat{\epsilon}{\ot}\i_C) (m{\ot}\i_{C{\ot}C}) (\i_B{\ot}F{\ot}\i_{C}) (\P{\ot}\P) \D\\
=&(m{\ot}\hat{\epsilon}{\ot}\i_C) (\i_B{\ot}F{\ot}\i_{C}) (\P{\ot}\P) \D\\
=&(m{\ot}\i_C) (\i_B{\ot}\hat{\epsilon}{\ot}\i_{B}{\ot}\i_C) (\P{\ot}\P) \D\\
=&(m{\ot}\i_C) ([(\i_B\ot\hat{\epsilon})\P]\ot\P) \D\\
=&(m{\ot}\i_C) (\omega\ot\P) \D\\
=&(m{\ot}\i_C) (\i_B\ot\P) (\omega\ot\i_A) \D
\end{split}
\end{equation*}
Since $\epsilon$ is a left counit for $(A,\D)$, we have
$$(\omega\ot\i_A)\D(a)=1_B\ot a,$$
for every $a\in A$. This implies that
\begin{equation*}
\begin{split}
(m{\ot}\i_C) (\i_B\ot\P) (\omega\ot\i_A) \D(a)&=(m{\ot}\i_C) (\i_B\ot\P)(1_B\ot a)\\
&=\P(a),
\end{split}
\end{equation*}
for every $a$ in $A$. This completes the proof.
\end{proof}
Analogous of Theorem \ref{t2} is satisfied for quantum groups that have right and (two-sided) counits.
Some natural questions about the structure of the compact quantum semigroup $(C,\G)$ arise:
\begin{question}
Let $(A,\D)$ and $(C,\G)$ be as in Theorem \ref{t1}.
\begin{enumerate}
\item [(i)] Suppose that $(A,\D)$ has one of the left or weak left cancelation properties.
Dose this hold for $(C,\G)$? In particular:
\item [(ii)] Suppose that $(A,\D)$ is a compact quantum group. Is $(C,\G)$ a compact quantum group?
\item [ (iii)] Are the converse of (i) and (ii) satisfied?
\end{enumerate}
\end{question}
We consider some parts of these questions for a simple example in the next section.
%%%%%%%%%%%%%%%%%%%%%%%%%%%%%%%%%%%%%%%%%%%%%%%%%%%%%%%%%%%%%%%%%%%%%%%%%%%%%%%%%%%%%%%%%%%%%%%%%%%%%%%%%%%%%%%%%%%%%%%%%%
\section{Some examples}
In this section, we consider a class of examples.
Let $A=\mathbb{C}^n$ be the C*-algebra of functions on the
commutative finite space $\{1,\cdots,n\}$, and let $(C,\P)$ be the quantum
family  of all maps from $\Q A$ to $\Q A$. A direct computation shows that $C$
is the universal C*-algebra generated by $n^2$ elements $\{c_{ij}: 1\leq i,j\leq n\}$ that satisfy the
relations
\begin{enumerate}
\item [(1)] $c_{ij}^2=c_{ij}=c_{ij}^*$ for every $i,j=1,\cdots,n$,
\item [(2)] $\sum_{j=1}^nc_{ij}=1$ for every $i=1,\cdots,n$, and
\item [(3)] $c_{ij}c_{ik}=0$ for every $i,k,j=1,\cdots,n$.
\end{enumerate}
Also, $\P:A\to A\ot C$ is defined by
$\P(e_k)=\sum_{i=1}^ne_i\ot c_{ik}$,
where $e_1,\cdots,e_n$ is the standard basis for $A$. Suppose that
$$\xi:\{1,\cdots,n\}\times\{1,\cdots,n\}\to\{1,\cdots,n\}$$
is a semigroup multiplication. Then $\xi$ induces a comultiplication $\D: A\to A\ot A$
$$\D(e_k)=\sum_{r,s=1}^n\D_k^{rs}e_r\ot e_s,$$
defined by $\D_k^{rs}=\delta_{k\xi(r,s)}$, where $\delta$ is the Kronecker delta.
We compute the comultiplication $\G:C\to C\ot C$, induced by $\D$ as in
Theorem \ref{t1}. We have
\begin{equation*}
\begin{split}
(\P\ot\P)\D(e_k)&=(\P\ot\P)(\sum_{r,s=1}^n\D_k^{rs}e_r\ot e_s)=\sum_{r,s=1}^n\D_k^{rs}\P(e_r)\ot \P(e_s)\\
&=\sum_{r,s=1}^n\sum_{j=1}^n\sum_{i=1}^n\D_k^{rs}e_j\ot c_{jr}\ot e_i\ot c_{is},
\end{split}
\end{equation*}
and therefore
\begin{equation*}
\begin{split}
&(m{\ot} \i_{C{\ot} C}) (\i_B{\ot} F{\ot} \i_C) (\P{\ot}\P) \D(e_k)\\
=&(m{\ot} \i_{C{\ot} C})(\sum_{j=1}^n\sum_{i=1}^n\sum_{r,s=1}^n\D_k^{rs}e_j\ot e_i\ot c_{jr}\ot  c_{is})\\
=&\sum_{l=1}^n\sum_{r,s=1}^n\D_k^{rs}e_l\ot c_{lr}\ot  c_{ls}=\sum_{l=1}^ne_l\ot(\sum_{r,s=1}^n\D_k^{rs}c_{lr}\ot  c_{ls}).
\end{split}
\end{equation*}
This equals to
$(\i_A\ot\G)\P(e_k)=(\i_A\ot\G)\sum_{i=1}^ne_i\ot c_{ik}=\sum_{l=1}^ne_l\ot\G(c_{lk})$.
Thus $\G$ is defined by
$$\G(c_{lk})=\sum_{r,s=1}^n\D_k^{rs}c_{lr}\ot  c_{ls}.$$
We now consider the special case $n=2$, in more details.
There are only five semigroup structures (up to isomorphism) on the set $\{1,2\}$:
\begin{enumerate}
\item [$\xi_1:$] $11=1,\hspace{2mm}12=2,\hspace{2mm}21=2,\hspace{2mm}22=1$.
\item [$\xi_2:$] $11=1,\hspace{2mm}12=2,\hspace{2mm}21=2,\hspace{2mm}22=2$.
\item [$\xi_3:$] $11=1,\hspace{2mm}12=1,\hspace{2mm}21=1,\hspace{2mm}22=1$.
\item [$\xi_4:$] $11=1,\hspace{2mm}12=1,\hspace{2mm}21=2,\hspace{2mm}22=2$.
\item [$\xi_5:$] $11=1,\hspace{2mm}12=2,\hspace{2mm}21=1,\hspace{2mm}22=2$.
\end{enumerate}
For every semigroup $(\{1,2\},\xi_i)$, let $(C,\G_i)$ be the corresponding quantum semigroup, as above.
A simple computation shows that:

$\Bigg\{
        \begin{array}{cc}
          \G_1(c_{11})=c_{11}\ot c_{11}+c_{12}\ot c_{12} &  \G_1(c_{12})=c_{11}\ot c_{12}+c_{12}\ot c_{11}\\
          \G_1(c_{21})=c_{21}\ot c_{21}+c_{22}\ot c_{22} &  \G_1(c_{22})=c_{21}\ot c_{22}+c_{22}\ot c_{21}\\
        \end{array}
        $

$\Bigg\{
        \begin{array}{cc}
          \G_2(c_{11})=c_{11}\ot c_{11} &  \G_2(c_{12})=c_{12}\ot c_{12}+c_{11}\ot c_{12}+c_{12}\ot c_{11}\\
          \G_2(c_{21})=c_{21}\ot c_{21} &  \G_2(c_{22})=c_{22}\ot c_{22}+c_{21}\ot c_{22}+c_{22}\ot c_{21}\\
        \end{array}
        $

$\Bigg\{
        \begin{array}{cc}
          \G_3(c_{11})=1 & \G_3(c_{12})=0 \\
          \G_3(c_{21})=1 & \G_3(c_{22})=0 \\
        \end{array}
        $

$$\G_4(c)=c\ot1\hspace{5mm}(\forall c\in C)$$
$$\G_5(c)=1\ot c\hspace{5mm}(\forall c\in C)$$
Note that the semigroup structure $\xi_1$ is a group structure and, $\xi_4$ and $\xi_5$ has right and left cancelation
properties, respectively. Thus $(\mathbb{C}^2,\D_1)$ is a compact quantum group and, $(\mathbb{C}^2,\D_4)$ and $(\mathbb{C}^2,\D_5)$
are compact quantum semigroups with right and left cancelation properties, respectively. From above computations, it is clear that
the compact quantum semigroups $(C,\G_4)$ and $(C,\G_5)$ have right and left cancelation properties, respectively. Now, we show that
$(C,\G_1)$ is also a compact quantum group: The unital C*-algebra $C$ is generated by the two unitary elements $x=c_{11}-c_{12}$ and
$y=c_{21}-c_{22}$ (see the following Remark for more details). A simple computation shows that
$$\G_1(x)=x\ot x\hspace{5mm}\text{and}\hspace{5mm}\G_1(y)=y\ot y.$$
This easily implies that $(C,\G_1)$ has left and right cancelation properties, and therefore $(C,\G_1)$
is a compact quantum group.
\begin{remark}
\begin{enumerate}
\item [(1)] The algebra $A=\mathbb{C}^2$ is the universal C*-algebra
generated by a unitary self-adjoint element, say $(1,-1)$.
It follows from the proof of Theorem 3.3 of \cite{S1}, that $C$ becomes the universal
C*-algebra generated by two unitary self-adjoint elements.
A model for $C$, is the C*-algebra of all continuous maps from  closed unit interval to  $2\times2$ matrix algebra,
which get diagonal matrices at the endpoints of the interval, equivalently
$$C=\Bigg\{
      \left(
        \begin{array}{cc}
          f_{11} & f_{12} \\
          f_{21} & f_{22} \\
        \end{array}
      \right)
:f_{ij}\in\mathcal{C}[0,1],f_{12}(0)=f_{12}(1)=f_{21}(0)=f_{21}(1)=0\Bigg\},$$
with unitary self-adjoint generators
$$x=
    \left(
      \begin{array}{cc}
        \cos(\pi t) & \sin(\pi t) \\
        \sin(\pi t) & -\cos(\pi t) \\
      \end{array}
    \right)\hspace{5mm}\text{and}\hspace{5mm}
y=
    \left(
      \begin{array}{cc}
        -\cos(\pi t) & \sin(\pi t) \\
        \sin(\pi t) & \cos(\pi t) \\
      \end{array}
    \right).
$$
In this representation of $C$, the generators $c_{ij}$'s become:
$c_{11}=\frac{1+x}{2}$, $c_{12}=\frac{1-x}{2}$, $c_{21}=\frac{1+y}{2}$ and $c_{22}=\frac{1-y}{2}$.
Also, the homomorphism $\P:\mathbb{C}^2\to \mathbb{C}^2\otimes C=C\oplus C$ is defined by $\P(1,-1)=(x,y)$.
This representation of the C*-algebra $C$, is one of the elementary examples of non-commutative spaces, see section
II.2.$\beta$ of \cite{C}.
\item [(2)]There is another quantum semigroup structure on quantum families of all maps from any
finite quantum space to itself introduced by So{\l}tan \cite{S1}.
\end{enumerate}
\end{remark}
\textbf{Acknowledgement:}
The author is grateful to the referee for his/her valuable suggestions.
%%%%%%%%%%%%%%%%%%%%%%%%%%%%%%%%%%%%%%%%%%%%%%%%%%%%%%%%%%%%%%%%%%%%%%%%%%%%%%%%%%%%%%%%%%%%%%%%%%%%%%%%%%%%%%%%%%%%%%%%%%

%%%%%%%%%%%%%%%%%%%%%%%%%%%%%%%%%%%%%%%%%%%%%%%%%%%%%%%%%%%%%%%%%%%%%%%%%%%%%%%%%%%%%%%%%%%%%%%%%%%%%%%%%%%%%%%%%%%%%%%%%%
\end{document}